\documentclass[12pt,a4paper,amsmath,amssymb,graphicx,epsfig,subfigure]{article}
\usepackage{a4wide}
\usepackage{amsmath,amssymb, textcomp}
\usepackage{epsfig}
\usepackage{float}
\usepackage{comment}
\usepackage{textcomp}
\usepackage[utf8]{inputenc}
\usepackage{gensymb}
\usepackage{tikz}
\usepackage{comment}
\usepackage{multirow}
\usepackage{latexsym}
\usepackage{amssymb}
\usepackage{amsmath}
\usepackage{tikz}
\usetikzlibrary{calc}

\tikzstyle{vertex} = [circle, draw=black, fill=black, inner sep=0.5mm, outer sep=0mm]
\tikzstyle{invvertex} = [circle, draw=white, fill=white, inner sep=0.0mm, outer sep=0mm]
\usetikzlibrary{shapes.geometric, arrows, positioning}
\usepackage{geometry}
\tikzstyle{title} = [rectangle, rounded corners, minimum width = 5cm, minimum height = 1cm, text centered, draw = black, fill=red!30]
\tikzstyle{mathematician} = [rectangle, rounded corners, minimum width = 3cm, minimum height = 1cm, text centered, draw = black, fill=red!10]

\tikzstyle{arrow} = [thick, ->, >=stealth]

\usetikzlibrary{quotes, angles, calc}

\usepackage{float}
\usepackage{comment}
\usepackage{textcomp}
\usepackage{hhline}
\usepackage{draftcopy}
\draftcopyName{VIJAY A. SINGH}{95}

\newcommand{\BEA}{\begin{eqnarray}}
\newcommand{\EEA}{\end{eqnarray}}
\newcommand{\BEAs}{\begin{eqnarray*}}
\newcommand{\EEAs}{\end{eqnarray*}}

\begin{document}
\begin{center}
  \begin{large}\textbf{Aryabhata and the Construction of the First Trigonometric Table}
  \end{large}
\end{center}
\vskip 0.5cm
\begin{center}
\textbf{Vijay A. Singh}\footnote{emailid: physics.sutra@gmail.com}\\
\textbf{Visiting Professor, Physics Department, Centre for for Excellence in Basic Sciences, Kalina, Santa-Cruz East, Mumbai-400098, India}
\vskip 0.5cm
\textbf{Aneesh Kumar}\\
\textbf{Dhirubhai Ambani International School, BKC, Bandra East, Mumbai-400098, India}\\
%\textbf{Email: physics.sutra@gmail.com}\\ \vspace{0.1in}
%\textit{(Jan. 1 2023 updated April 1 2023)}\\,\\
\vskip 0.5cm
\end{center}
 Few among us  would know that the  first mention of the  sine and the
 enumeration  of  the   first  sine  table  are  to   be  credited  to
 Aryabhata. The method to generate  this relies on the sine difference
 formula which is  derived using ingenious arguments  based on similar
 triangles.  We describe  how  this was  done. In  order  to make  our
 presentation pedagogical we take a unit circle and radians instead of
 the  (now)  archaic  notation  in the  \textit{Aryabhatiya}  and  its
 commentators. We suggest a couple of exercise problems and invite the
 enterprising student to  try their hands. We also point  out that his
 sine and  the second  sine difference identities  are related  to the
 finite  difference  calculus  we   now  routinely  use  to  calculate
 derivatives  and  second  derivatives.   An  understanding  of  these
 trigonometric  identities  and preparation  of  the  sine table  will
 enable a student to get an  appreciation of the path breaking work of
 Aryabhata.

\section{Introduction}

We routinely look  up a table of  sines or punch on  our calculator to
obtain the value of a trigonometric function for a given angle. Little
do we realize  that a debt of  gratitude is owed to  the fifth century
Indian savant Aryabhata who first showed the way. Aryabhata enumerated
the table of  sines for closely spaced angles. His  methods were based
on   general  trigonometric   identities   and   lend  themselves   to
extensions. The first  mention of the sine function is  to be found in
his (one  and only)  seminal work the  \textit{Aryabhatiya} (499\,CE).
Aryabhata  describes it  in poetic  terms  as the  half bow-string  or
\textit{Ardha-Jya}.  This  is illustrated  in Fig.   1.  The  arrow or
\textit{saar} is related  to the cosine function. This is  not the only
example of  poetry intruding into  his mathematics. To describe  the fact,
heretical  and revolutionary  for  those  times and  for  a long  time
afterwards,  that the  earth is  rotating and  the sun  is stationary,
Aryabahta evokes  the tranquil  metaphor of a  boat floating  down the
river  and the  stationary land  mass  which  seems to  move
backwards. The  \textit{Aryabhatiya} is  in verse  form and  there are
over a  100 of them.  They  are written fully respecting  the norms of
grammar and metre. Aryabhata was a poet as well.

The  \textit{Aryabhatiya} consists  of 121  cryptic verses,  dense and
laden with  meaning \cite{kern,shukla}.   The work  is divided  into 4
parts  or \textit{padas}:  the  \textit{Gitikapada}  (13 verses),  the
mathematics     or     \textit{Ganitapada}    (33     verses),     the
\textit{Kalkriyapada}    (25   verses)    and    the   astronomy    or
\textit{Golapada} (50 verses).  The  astronomy is better known.  There
are two  verses in the mathematics  \textit{Ganitapada} describing the
solution of  the linear  Diophantine equation.  This has  received due
recognition.   Our focus  here  is  on the  trigonometry  part in  the
\textit{Ganitapada} which in our opinion has suffered neglect and is a
pioneering achievement of this savant.

Virtually  every  major Indian  mathematician  has  commented on  the
\textit{Aryabhatiya}  (499 CE).   Often it  is  in terms  of a  formal
\textit{Bhasya}  (Commentary). Table  I lists  some of  these. Notable
among them is  the voluminous work (\textit{Maha Bhashya})  of the 15th
century mathematician Nilakantha Somaiyaji (1444 CE - 1544 CE) who was
part of  the Kerala school  which, beginning  with Madhava (1350  CE -
1420  CE),  founded  the  calculus of  trigonometric  functions.   Our
presentation relies on  Somaiyaji's commentary \cite{soma} and the  works of Kripa
Shankar  Shukla  and   K.  V.  Sarma  \cite{shukla}  as   well  as  of
P. P. Divakaran \cite{ppd}.

In  this article  we  describe the  trigonometric  identities used  by
Aryabhata  to obtain  the table  of  sines.  This  entails taking  the
difference of  the sine of two  closely spaced angles and  then taking
the second  sine difference.  We  next show that these  identities are
the  same as  the  finite  difference calculus  one  uses nowadays  to
numerically obtain the first derivative  and second derivative of sine
functions. We follow this up with a brief discussion.
\begin{comment}
We argue that this method  is general and
demonstrate it  by extending it  to obtain the corresponding  table of
cosines. We also
explore  the  expansions for  the  trigonometric  functions which  the
Aryabhata  identities lead  to. These  are different  from the  Taylor
series  we are  now familiar  with.  
\end{comment}

The Indian mathematical  tradition is largely word  based. Results are
mentioned and derivations are omitted.   The \textit{Aryabhatiya} (499 CE) with
a  little  over  a  100  cryptic,  super-compressed  verses  of  dense
mathematics is a  prime example.  Our approach is  pedagogical and one
which will help the student  and teacher to appreciate this pioneering
work.   Hence we  shall take  some  liberties and  describe the  great
master's work in  terms of unit circle and radians  instead of degrees
and minutes. A couple of exercises in the end will help one get a more
hands on understanding.

%\newpage

\begin{table}[H]
    \centering
    \begin{tabular}{|c|c|c|c|}
        \hline Bhaskara I & 629 CE  & Sanskrit & Valabhi, Gujarat \\
        
        \hline Suryadeva Yajvan & born 1191 CE& Sanskrit & Gangaikonda-Colapuram\\
        
        \hline Parameshvara & 1400 CE  & Sanskrit & Allathiyur, Kerala \\
        
        \hline Nilakantha Somayaji & 1500 CE & Sanskrit & Trikandiyur, Kerala \\
        
        \hline Kondadarma & unknown & Telegu & Andhra \\

%          \hline Al-Khwarizmi & 820 CE & Arabic & Baghdad \\
        
        \hline Abul Hasan Ahwazi & 800 CE & Arabic & Ahwaz, Iran \\
        
        \hline
    \end{tabular}
    \caption{A host of eminent mathematicians have commented on the \textit{Aryabhatiya} written CE 499. Some, like Brahmagupta (600 CE) or Bhaskara II (1100 CE) have not written a specific commentary but  dwelt extensively on it. The above is an abbreviated list of specific commentaries and the dates are approximate. Our main source for this is the work of K. S. Shukla and K. V. Sarma \cite{shukla} which cites around 20 commentaries.  \label{table:1}}   
\end{table} 

\section{The \textit{Ardha-Jya} or Sine Function}
 
 As mentioned earlier the \textit{Aryabhatiya} has some 121 verses out
 of   which   33   verses    belong   to   the   mathematics   section
 (\textit{Ganitapada}). Aryabhata  works with,  for the first  time in
 the history of  mathematics, the sine function. It  is the half-chord
 $AP$  of  the   unit  circle  in  Fig.  1.    \BEAs  sin(\theta)  &=&
 \dfrac{AP}{OA} \\ &=&  AP \,\,\,\,\,\,\,\, (OA = 1)  \EEAs The circle
 maybe large  or small; correspondingly  $AP$ and $OA$ maybe  large or
 small,   but  the   l.h.s.  is   a  function   of  $\theta$   and  is
 \textbf{invariant}. Further,  all metrical properties related  to the
 circle  can   be  derived  using  trigonometric   functions  and  the
 Pythagorean  theorem (also  described as  the Baudhayana  or Diagonal
 theorem \cite{ppd}). For example  the
 geometric property  of a triangle can  be related to the  arcs of the
 circumscribing  circle using  the sine  and cosine  functions. Or  the
 diagonals  of the  inscribed  quadrilateral can  be  related to  its
 sides. (If a recent proof  of the this Pythagorean theorem using the
 the law of sines holds up to scrutiny  then \textit{all} metrical
 properties  of  a  circle  can  by  obtained  by  trigonometry  alone
 \cite{shirali}.)  By emphasizing  the role of this  half-chord Aryabhata endowed
 circle geometry with metrical properties.  This alone may qualify him
 as the founder of trigonometry. But he did more.

\begin{figure}[h!]
  \centering
  \includegraphics[scale=0.2]{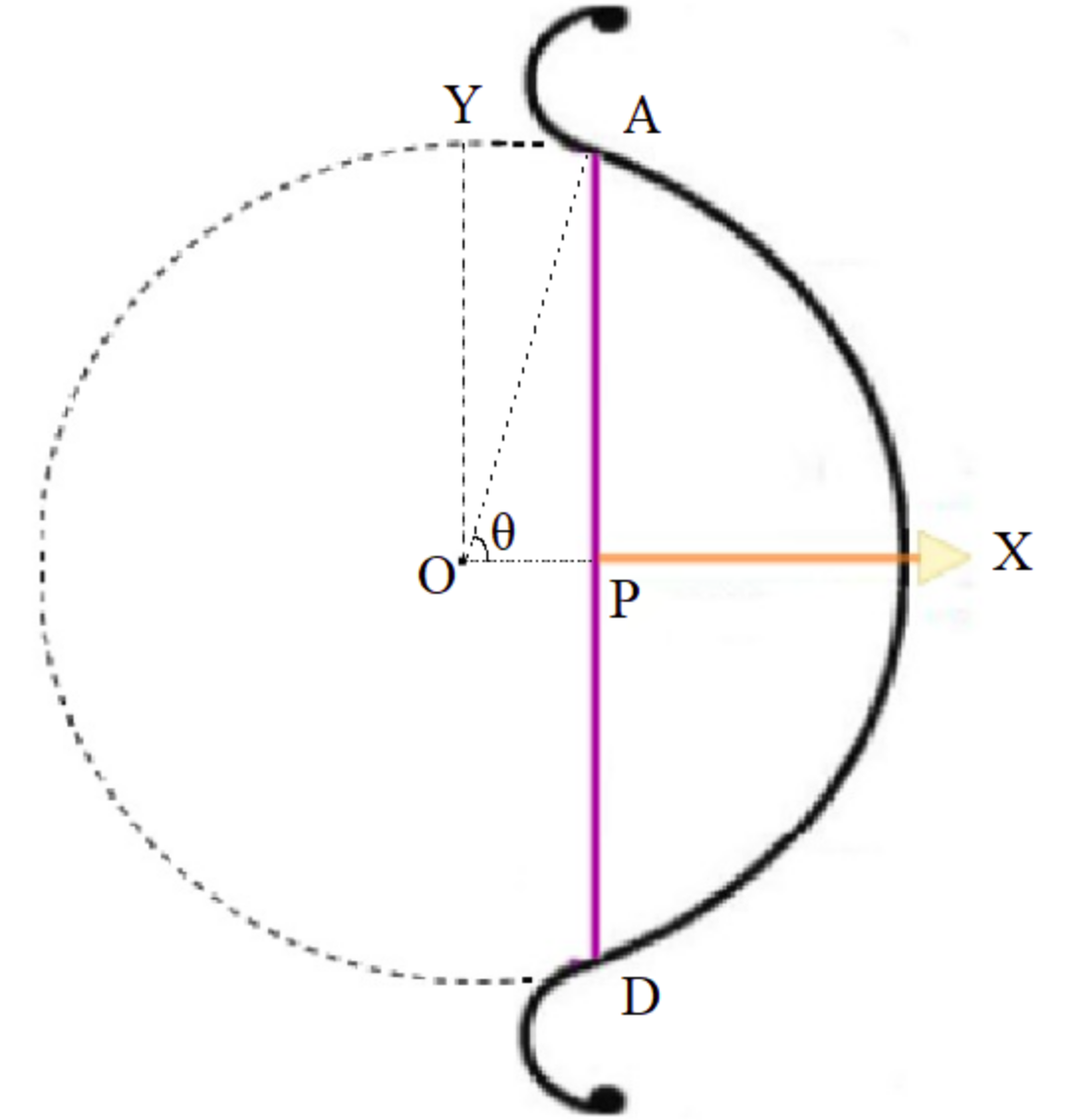}
\caption{The bow superposed on the unit circle. Half the  bow string or half chord $AP$ is $sin(\theta)$ as defined by Aryabhata. $OP$ is $cos(\theta)$ while $PX$ = 1-$cos(\theta)$ is called the \textit{saar}. See text for comments.}  
\end{figure}  

It may also help to note that the length of the half-chord $AP$ is very close to the length of the arc $AX$ only when the angle is small (e.g. sin($\theta$) $\approx \theta$ if $\theta$ is small and in radians). This was known to Aryabhata in all likelihood by inspection. Similarly the sine of 90$^{\degree}$ is 1, since then the half-chord is the same as the radius. In the 9th verse of the  \textit{Ganitapada} he uses the property of an equilateral triangle and  obtains the sine of 30$^{\degree}$ as 1/2.

We pause  to note  that Aryabhata  also states the  value of  $\pi$ as
62832/20000 in the  10th verse. This value is 3.416  and he is careful
enough to state  that this is proximate  (\textit{Asanna}) which means
we can obtain  better and truer values for $\pi$  presumably with more
effort\footnote{The  word  \textit{Asanna}  or   proximate  is  to  be
distinguished  from \textit{Sthula}  which is  approximate or  roughly
equal}.

\section{The Difference Formula for Sine and Cosine}

The 12th verse of the \textit{Ganitapada} plays a central role in the tabulation of the sine function. It is cryptic and to unravel its meaning we first need to obtain the difference formula for the sine. The presentation below relies on a number of sources: (i) The commentary of Nilakantha Somaiyaji \cite{soma}; (ii) the treatment of Shukla and Sarma \cite{shukla}; (iii) and for the sake of ease of understanding we follow Divarkaran \cite{ppd}  and take a unit circle as opposed to a radius of 3438 by earlier workers\footnote{One radian is 3438 minutes and we remind the reader that 2\,$\pi$ radians is 180 degrees and 1 degree is 60 minutes.}.

Figure 2 depicts a quadrant of the unit circle where $OX$ = $OY$ = 1. The arcs $XA$, $XB$ and $XC$ trace angles $\theta$, $\theta + \phi$ and $\theta - \phi$ respectively. The half-chords $AP$, $BQ$ and $CR$ are the corresponding sine functions. We drop a perpendicular $CS$ from the circumference onto the half-chord $BQ$ as shown. According to his commentator Nilakantha Somaiyaji \cite{soma}, Aryabhata obtained the  relationship  between the difference in the trigonometric functions by demonstrating that the two triangles $BSC$ and $OPA$ are similar and  then exploiting this ingenuously. We trace his line of reasoning in the Appendix where we derive the following relations

\BEA
sin(\theta + \phi) - sin(\theta - \phi) & = & 2 sin(\phi) cos(\theta) \label{sindiff1} 
\EEA

The difference in the sines is proportional to the cosine of the mean angle.

\BEA
cos(\theta + \phi) - cos(\theta - \phi) & = & - 2 sin(\phi) sin(\theta) \label{cosdiff1}
\EEA

The difference in the cosines is proportional to the (negative) of the sine of the mean angle.

%\begin{comment}
\begin{figure}
%  \begin{center}
     \begin{tikzpicture}[>=latex]
%  \begin{tikzpicture}[scale=0.6]
 % Theta and Phi Values
\pgfmathsetmacro{\ThetaVal}{50} % degrees
\pgfmathsetmacro{\PhiVal}{10} % degrees
% end Theta and Phi Values

% 'Theta - Phi' and 'Theta + Phi' label text angles (degrees)
\pgfmathsetmacro{\labeltextangle}{-65}
% end 'Theta - Phi' and 'Theta + Phi' label text angles (degrees)

\draw [draw=white](10,10) grid (20,20);
% drawing the figure -----------------------------------------------------------------------------------------------
% variable declaration
\pgfmathsetmacro{\BX}{{10 + 10*cos(\ThetaVal+\PhiVal)}}%{14.995}
\pgfmathsetmacro{\BY}{{10 + 10*sin(\ThetaVal+\PhiVal)}}%{18.663}
\pgfmathsetmacro{\CX}{{10 + 10*cos(\ThetaVal-\PhiVal)}}
\pgfmathsetmacro{\CY}{{10 + 10*sin(\ThetaVal-\PhiVal)}}
\pgfmathsetmacro{\AX}{{10 + 10*cos(\ThetaVal)}}
\pgfmathsetmacro{\AY}{{10 + 10*sin(\ThetaVal)}}

\coordinate (O) at (10,10);
\coordinate (X) at (20,10);
\coordinate (Y) at (10,20);
\coordinate (A) at (\AX,\AY);
\coordinate (B) at (\BX,\BY);
\coordinate (C) at (\CX,\CY);
\coordinate (Q) at (\BX,10);
\coordinate (P) at (\AX,10);
\coordinate (R) at (\CX,10);
\coordinate (S) at (\BX,\CY);

\pgfmathsetmacro{\startX}{0}
\pgfmathsetmacro{\endX}{90}
\pgfmathsetmacro{\radiusX}{10}
% end variable declaration
\draw (X) arc (\startX:\endX:\radiusX);
\draw (O) -- (X);
\draw (O) -- (Y);
\draw [style=dashed](O) -- (A);
\draw [style=dashed](O) -- (B);
\draw [style=dashed](O) -- (C);
\draw [style=dashed](B) -- (C);
\draw (Q) -- (B);
\draw (R) -- (C);
\draw (P) -- (A);
\draw (S) -- (C);
% labelling the points
\node [vertex, label=below:\Large \textbf{$O$}, name=O] at (O) {};
\node [vertex, label=right:\Large \textbf{$X$}, name=X] at (X) {};
\node [vertex, label=above:\Large \textbf{$Y$}, name=Y] at (Y) {};
\node [vertex, label=above right:\Large \textbf{$A$}, name=A] at (A) {};
\node [vertex, label=above right:\Large \textbf{$B$}, name=B] at (B) {};
\node [vertex, label=above right:\Large \textbf{$C$}, name=C] at (C) {};
\node [vertex, label=below:\Large \textbf{$Q$}, name=Q] at (Q) {};
\node [vertex, label=below:\Large \textbf{$P$}, name=P] at (P) {};
\node [vertex, label=below:\Large \textbf{$R$}, name=R] at (R) {};
\node [vertex, label=below left:\Large \textbf{$S$}, name=S] at (S) {};
% end of drawing the figure -----------------------------------------------------------------------------------------------

% drawing the angle and perpendicular symbols -----------------------------------------------------------------------------------------------
% variable declaration
\pgfmathsetmacro{\perpsize}{0.2}
% end of variable declaration
\draw (S) rectangle ++(\perpsize,\perpsize) ; % drawing the square denoting the perpendicular symbol.
\draw pic[draw, angle radius=2cm, angle eccentricity=1.5, "\Large\textbf{$\theta$}", <->, pic text options={shift={(-0.5,-0.5)}}] {angle=X--O--A};
\draw pic[draw, angle radius=2.9cm, angle eccentricity=1.5, "\Large\textbf{$\phi$}", <->, pic text options={shift={(-0.6,-0.8)}}] {angle=A--O--B};
\draw pic[draw, angle radius=3.1cm, angle eccentricity=1.5, "\Large\textbf{$\phi$}", <->, pic text options={shift={(-0.8,-0.8)}}] {angle=C--O--A};
% end of drawing the angle and perpendicular symbols -----------------------------------------------------------------------------------------------

% Arcs and labels for theta +/- phi -----------------------------------------------------------------------------------------------
\draw pic[draw, angle radius=12cm, angle eccentricity=1, "\Large\textbf{$\theta - \phi$}", |<->|, pic text options={shift={(-0.5,-0.5)}, rotate={\labeltextangle}}] {angle=X--O--C};
\draw pic[draw, angle radius=14cm, angle eccentricity=1, "\Large\textbf{$\theta + \phi$}", |<->|, pic text options={shift={(-0.5,-0.5)}, rotate={\labeltextangle}}] {angle=X--O--B};
% end of arcs and labels for theta +/- phi  -----------------------------------------------------------------------------------------------

    \end{tikzpicture}
    \caption{Derivation of the sine difference relation. The figure depicts the quadrant of a unit  circle of  radii $OX$ = $OY$ = 1. The half-chords $AP$, $BQ$ and $CR$ are $sin(\theta)$,   $sin(\theta + \phi)$ and  $sin(\theta - \phi)$ respectively. It is worth noting that (later) we shall take $\phi$ to be a small angle.}
%\end{center}
  \end{figure}
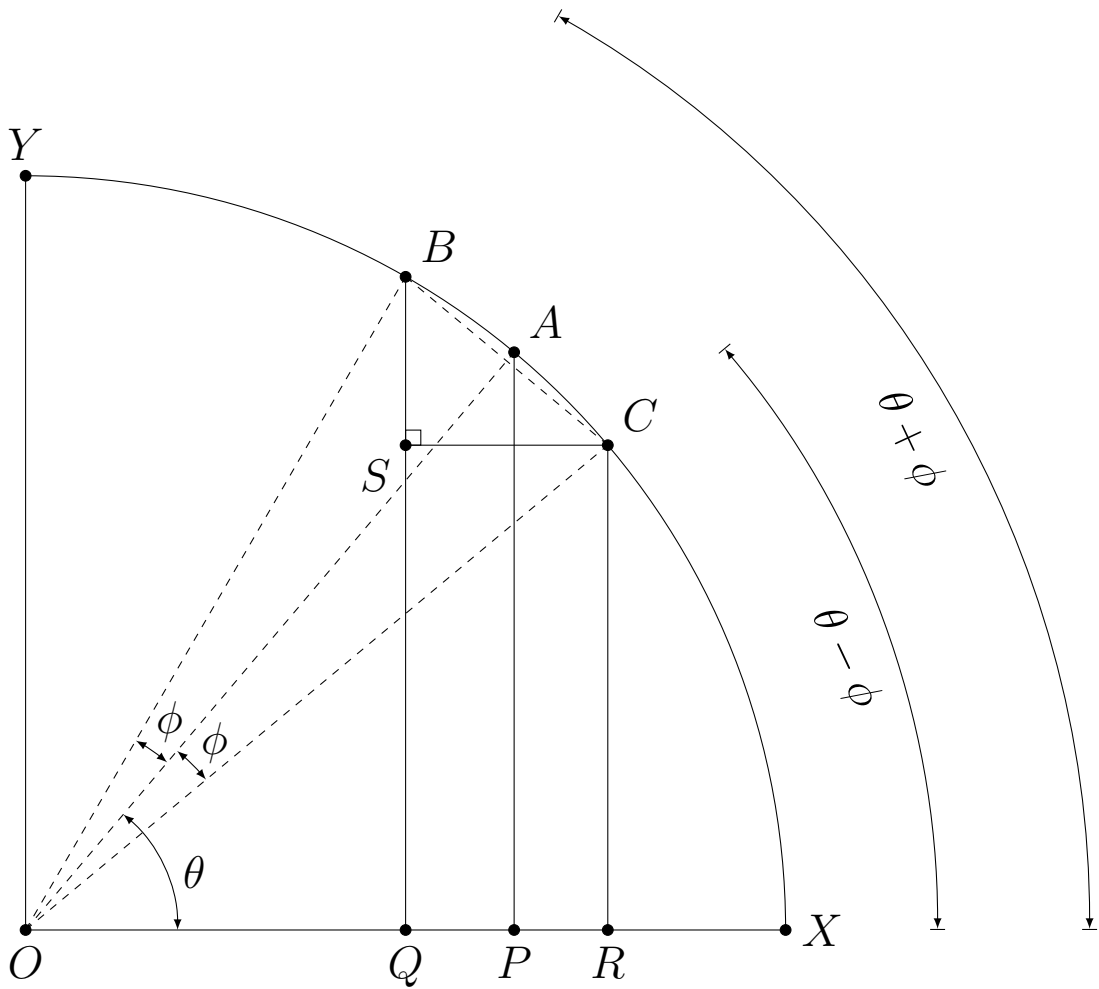   
%\end{comment}

\begin{comment}
\begin{figure}[h!]
  \centering
  \includegraphics[scale=0.4]{quadrant_unit_x.pdf}
%  \includegraphics{quadrant_unut_x.pdf}
\caption{Derivation of the sine difference relation. The figure depicts the quadrant of  unit circle of radii $OX$ = $OY$ = 1. The half-chords $AP$, $BQ$ and $CR$ are $sin(\theta)$,   $sin(\theta + \phi)$ and  $sin(\theta - \phi)$ respectively. It is worth noting that later we shall take $\phi$ to be a small angle.}  
\end{figure}
\end{comment}

%\newpage
\section{The Sine Table}

Aryabhata  obtained the values of the sines at fixed angles between  0 and $\pi /2$ thus  generating the sine table for $\pi/48$ = 3.75 degrees, 7.5 degrees up to 90 degrees. This table is stated in verse 12 of the first chapter, the \textit{Gitikapada}. The table has been used by Indian astronomers (and astrologers) in some form or another since 499 CE up to the present. We shall see how the table was generated.

Let us take $\phi = \epsilon/2$  where $\epsilon$ is small. We take $\theta = (n-1/2) \epsilon$ where $n$ is a positive integer from 1 to $N$.  To fix our ideas $\epsilon = \pi$/48 = 3.75$^0$ = 225' and $N$ = 24. We re-state the sine and cosine difference formulae (Eqs.\,(\ref{sindiff1}) and (\ref{cosdiff1})) from the previous section 
\BEA
\delta s_{n} &=& s_{n} - s_{n-1} = 2 s_{1/2} c_{(n-1/2)} \label{finite1} \\
 \delta c_n &=& c_{n} - c_{n-1} = -2 s_{1/2} s_{(n-1/2)} \label{finite2} 
 \EEA
 where the symbol $s_n$ stands for $sin(n \epsilon)$, $c_n$ stands for $cos(n \epsilon)$ and $s_{1/2}$ for $sin(\epsilon/2)$. The above is a pair of coupled equations and it was Aryabhata's insight to take the second difference, namely
\BEA
\delta^2 s_n &=& \delta s_{n} - \delta s_{n-1}  = 2 s_{1/2} (c_{(n-1/2)}-c_{(n-3/2)})  \nonumber  \\
& =& - 4 s_{1/2}^2 s_{n-1} \mbox{\,\,\,\,\,\,\, on using Eq.\,(\ref{finite2}) }  \label{finite3}
 \EEA 
 Thus the second difference of the sines is proportional to the sine itself. The next step is to represent the r.h.s in terms of a recursion. We observe $s_n$ on the r.h.s. of Eq.(\ref{finite3}) may be written as $s_n = s_n - s_{n-1} + s_{n-1} - s_{n-2} + s_{n-2} - ...$ = $\delta s_{n-1} + \delta s_{n-2} + ...$  Thus
\BEA 
 \delta s_{n} - \delta s_{n-1} &=& - 4 s_{1/2}^2 \sum_{m=1} ^{n-1} \delta s_m \label{recur1}
 \EEA

 Thus we get a recursion relation where the second difference of the sines is expressed in terms of all previously obtained first sine differences. To initiate the recursion we need $\delta s_1$ which is $s_1 - s_0 = sin(\epsilon) - sin(0) \approx \epsilon$ since for small angles the half-chord and the corresponding arc are equal as stated in the previous section.

 Using the recursion relation of the sine, we can generate the celebrated sine table of Aryabhata, taking $\pi$ = 3.1416 and sin($\epsilon$) = $\epsilon$ = 0.0654 (= 225').

Table \ref{table:3} depicts some typical values of the sine function as well as the value of the sine multiplied by 3438 (the so called '\textit{Rsine}' of Aryabhata). We can see that this matches Aryabhata's celebrated sine table, up to $\pm$ 1 minute. For example, $\theta$ = $\pi$/6 gives 1719 minutes. For comparison, we also show the modern values of sine up to four decimal places. Note that Aryabhata takes angles up to $\pi$/2 and seemed aware of the fact that going further was unnecessary given the periodic nature of the sine function.

\begin{table}[H]
    \centering
    \begin{tabular}{|c|c|c|c|}
        \hline \hline $\theta$ & sin($\theta$) Aryabhata & sin($\theta$) (minutes) & sin($\theta$) Modern \\
\hline
\hline $\pi$/48 &   0.0654 & 225 & 0.0654 \\
\hline $2\pi$/48 & 0.1305  & 449 & 0.1305 \\
\hline $3\pi$/48 & 0.1951     &   671 & 0.1951 \\
\hline $4\pi$/48 & 0.2588     &   890 & 0.2588 \\
\hline $5\pi$/48 & 0.3214      &  1105 & 0.3214 \\
\hline $6\pi$/48 & 0.3827     &   1315 & 0.3827 \\
\hline $7\pi$/48 & 0.4423  &      1520 & 0.4423  \\
\hline $8\pi$/48 & 0.5000  &      1719 & 0.5000 \\
\hline $9\pi$/48 & 0.5556    &    1910 & 0.5556 \\
\hline $10\pi$/48 & 0.6088  &    2093 & 0.6088 \\
\hline $11\pi$/48 & 0.6594 &     2267 & 0.6593 \\
\hline $12\pi$/48 & 0.7072  &    2431 & 0.7071 \\
\hline $13\pi$/48 & 0.7519   &   2585 & 0.7518  \\
\hline $14\pi$/48 & 0.7935  &    2728 & 0.7934  \\
\hline $15\pi$/48 & 0.8316  &    2859 & 0.8315 \\
\hline $16\pi$/48 & 0.8662  &    2978 & 0.8660  \\
\hline $17\pi$/48 & 0.8971  &   3084 & 0.8969  \\
\hline $18\pi$/48 & 0.9241  &    3177 & 0.9239  \\
\hline $19\pi$/48 & 0.9472  &    3256 & 0.9469  \\
\hline $20\pi$/48 & 0.9662  &    3322 & 0.9659  \\
\hline $21\pi$/48 & 0.9812  &    3373 & 0.9808  \\
\hline $22\pi$/48 & 0.9919  &    3410 & 0.9914  \\
\hline $23\pi$/48 & 0.9983   &   3432 & 0.9979  \\
\hline $24\pi$/48 & 1.0005  &    3439 & 1.0000  \\
\hline \hline        
    \end{tabular}
    \caption{Table of values of sine  using the Aryabhatan method, taking episilon = $\pi$/48 (= 3.75$\degree$ = 225') and pi=3.1416 and comparison with the modern day values. In column 3 we also quote values in minutes as done in Verse 12,  \textit{Gitika} chapter of the \textit{Aryabhatiya} \cite{kern,shukla}.  \label{table:3}}  
\end{table} 

%\end{comment}

   \section{Finite Difference Calculus}

  Of greater relevance is the fact that the sine (or cosine) difference formulae foreshadow finite difference calculus, a popular numerical technique in this age of computation. Rewriting Eqs.\,(\ref{sindiff1}) and (\ref{cosdiff1}) with $\phi = \epsilon$,
  \BEA
  \dfrac{sin(\theta + \epsilon) - sin(\theta - \epsilon)}{2 sin(\epsilon)} & = & cos(\theta) \label{fdc1}  \\
    \dfrac{cos(\theta + \epsilon) - cos(\theta - \epsilon)}{2 sin(\epsilon)} & = & - sin(\theta) \label{fdc2} 
    \EEA
    
    Aryabhata took $\epsilon$ to be $\pi/48$. But he also stated that its value is \textit{yateshtani} or as per our wish (Verse 11, \textit{Ganitapada}). Some took it to be $\pi/96$ and others like Brahmagupta took it as $\pi/12$ or 15$^{\degree}$. If we take $\epsilon$ to be sufficiently small we have our classic formula for finite difference calculus. Noting that 2 $sin(\epsilon/2) \approx \epsilon$ we have the finite difference version of the derivative  of sine
      $$ \dfrac{\delta sin(\theta)}{\delta \theta} = cos(\theta)    $$
      and similarly for the cosine.
      $$ \dfrac{\delta cos(\theta)}{\delta \theta} = -sin(\theta)    $$
      
    Let us understand this with an example. We know that sine(37$^0$) is close to 0.6 and sine(30$^0$) is 0.5. The difference in angle is 7$^0$ which in radians is 0.122. Thus the derivative of sine of the median angle 33.5$^0$ from Eq.\,(\ref{fdc1})  is
 \BEAs
 \delta sin(\theta) /\delta \theta  &=& (0.6 - 0.5)/0.122 = .82
 \EEAs
 Looking up the sine table or the calculator yields cos(33.5) = 0.83. Similarly Eqs.\,(\ref{finite3})  yields the second derivatives namely
 \BEAs
 \delta^2 sin(\theta) / \delta^2 \theta 
      &  \approx & - sin(\theta)  \\
  \delta^2 cos(\theta) / \delta^2 \theta & \approx & - cos(\theta)
 \EEAs
 The above are now called central difference approximations to the derivative and the second derivative. Aryabhata does not mention the term  finite difference calculus (let alone calculus). But similar methods  are now  used to numerically solve our differential equations. A student can readily recognize the above as a standard solution of the classical simple harmonic oscillator. Note also that  Newton's II Law and the famous Schrodinger equation of quantum mechanics are both second order differential equations.

\section{Discussion}

One can discern a continuity in Indian mathematics, howsoever tenuous,  from pre-Vedic times ($<$ 1000 BCE) up and until 1800s. A striking example is the influence of the  \textit{Aryabhatiya} (499 CE) on major Indian mathematicians who followed him including Madhava (1350 CE) who founded Calculus \cite{ppd}; as also the influence on Aryabhata of the mathematics which preceded him \cite{datta,amma}.

Aryabhata describes the sine function  ``poetically'' as the \textit{Ardha-Jya} or half bow (see Fig. 1). To reiterate, he  seemed aware that the sine of (i) zero is zero; (ii) small angle is itself, or the small arc is almost equal to the half chord; (iii) 30 degrees is 1/2 (\textit{Ganitapada} verse 9); (iv) and 90 degrees is unity. Further, that the sine function is periodic so he  prudently  stops the enumeration of the sine for angles greater than 90 degrees. Then, in a remarkably insightful way, he laid down the recursion relation for sine differences which enables one to generate the sine table. It is this work, more than his solution to the linear Diophantine equation (verses 31 and 32 of the \textit{Ganitapada}) which establishes him as a genius and one of the brightest stars in the firmament of world mathematics.    

The sine table can also be generated using the half angle formula. This was demonstrated in the \textit{Panchasiddhantika} a text written barely 50 years after the appearance of the \textit{Aryabhatiya} \cite{vara}. As pointed out, a feature of the Aryabhata's difference relation is how contemporaneous it is. It can be easily recognized as finite difference calculus. It also led to the development of the calculus of trigonometric functions by the Madhava (1350 CE) and his disciples along the banks of the Nila river in Kerala. This school is variously called the Nila \cite{ppd} and even as the Aryabhata school \cite{amma}. Another aspect to note is that Bhaskara II (1100 CE) used the canonical 2 $\pi$ /96 division of the great circle to carry out discrete integration and obtain the (correct) expressions for the surface area and volume of the sphere. Jyesthdeva of the Nila (or Aryabhta) school in his work \textit{Yuktibhasa}  derived the same results using calculus (circa 1500 CE). Aryabhata can legitimately be called the founder of trigonometry.   

To sum up, the \textit{Arybhatiya} exercised a tremendous influence over Indian mathematicians and, for over a thousand years. For a book with just over a 100 pithy verses, its legacy remains unparalleled in the scientific world. We hope that our article will give our young audience an introduction to his work and will serve as an inspiration.

\paragraph{Acknowledgement:} One of the  authors (VAS)  would place on record the many useful discussions he had with Prof. P. P. Divakaran.

\appendix

\section*{Appendix: The Derivation of the sine and cosine difference formula}

As stated in the text Figure 2 depicts a quadrant of the unit circle where $OX$ = $OY$ = 1. The arcs $XA$, $XB$ and $XC$ trace angles $\theta$, $\theta + \phi$ and $\theta - \phi$ respectively. The half-chords $AP$, $BQ$ and $CR$ are the corresponding sine functions. We drop a perpendicular $CS$ from the circumference onto the half-chord $BQ$ as shown.

\begin{comment}The triangles $BCS$ and $OAP$ have sides which are mutually perpendicular; $CS$ is perpendicular to $BQ$, $BS$ is perpendicular to $OP$ and $BC$ is perpendicular to $OA$.
\end{comment}
  We  show that  the two  triangles $BSC$  and $OPA$  are similar.  By
  construction   $\angle{BSC}$   and   $\angle{OPA}$   are   each   90
  degrees. Note  $OB$ = $OC$  = 1  (unit radius) and  hence $\triangle
  OBC$ is isosceles.  This implies that $\angle{OBC}$  = 90$^{\circ} -
  \phi$.  Also  $\angle{OBQ}$ =  90$^{\circ} -  \phi -  \theta$. Hence
  $\angle{SBC}$=$\angle{OBC}$-$\angle{OBQ}$=$\theta$.       Therefore,
  $\angle{SBC}$=$\angle{POA}$=$\theta$.    This     establishes    the
  similarity of the two triangles by the angle-angle test.

\BEAs
\dfrac{BS}{OP} & = & \dfrac{BC}{OA}
\EEAs

Now $OA$ = 1 (unit circle), $OP$ = $cos(\theta)$ and $BC$ = 2 $sin(\phi)$. By inspection $BS$ = $BQ$ - $CR$ =  $sin(\theta + \phi)$ - $sin(\theta - \phi)$. This yields the sine difference formula (Eq.\,(\ref{sindiff1})

\BEAs
sin(\theta + \phi) - sin(\theta - \phi) & = & 2 sin(\phi) cos(\theta)  
\EEAs

Thus the difference in the sines is proportional to the cosine of the mean angle. We can also obtain the cosine difference formula (Eq.\,(\ref{sindiff1}) by noting that

\BEAs
\dfrac{CS}{AP} & = & \dfrac{BC}{OA}
\EEAs

Note $AP$ = $sin(\theta)$ and $CS$ = $OR$ - $OQ$ =  $cos(\theta - \phi)$ - $cos(\theta + \phi)$. Hence

\BEAs
cos(\theta + \phi) - cos(\theta - \phi) & = & - 2 sin(\phi) sin(\theta) 
\EEAs

The difference in the cosines is proportional to the (negative) of the sine of the mean angle. We pause to note that \textit{prima facie} the two triangles we considered appear unrelated. A hallmark of Indian mathematics is strong geometric intuition and this dates back to the \textit{Sulbasutra circa} 800 BCE. Another is the reliance on the ``rule of three'' (\textit{trirasikam}) - here we employ a simple version of it namely, if $a/b = c$ then $a = b \times c$.

\vskip 1.0cm
\begin{center}
  \textbf{EXERCISES}
\end{center}
\begin{enumerate}

\item We can generate the sine table as per Aryabhata's suggestion but not exactly using the same value for $\epsilon$ he used. We choose $\epsilon = \pi/80 \approx 0.0393$ which is the same as 2.25$^0$. We take $sin(\epsilon) \approx \epsilon$. If you have a simple calculator generate all values of sine from 2.25 to 18 degrees in equal steps using Eq.\,(\ref{recur1}). Alternatively if you have a programmable calculator or a computer generate all values of sine from 2.25 to 90 degrees. Compare with the results your calculator will otherwise yield.      

\item In the last section reference is made of the text  \textit{Panchasiddhantika} wherein the half angle forumla is mentioned, viz.
  $$ cos(2 \theta) = 1 - 2 sin^2(\theta)    $$
How would you (i) derive this by a geometrical construction; (i) employ this to generate the sine table?  
  
\end{enumerate}

%\newpage

%\begin{center}
%  \textbf{FEW REFERENCES}
%\end{center}
%These references are in English. \\
%\textbf{Primary}
%\begin{enumerate}

%\end{enumerate}

%\textbf{Secondary:}
%\begin{enumerate}

\vskip 1.0cm

\textit{Prof. Vijay A. Singh  has been faculty at IIT Kanpur (1984-2014) and HBCSE, Tata Institute for Fundamental Research (2005-2015) where he was the National Coordinator of  both the Science Olympiads and the National Initiative on Undergraduate Science for a decade. He is a Fellow, National Academy of Sciences, India and was President, the Indian Association of Physics Teachers (2019-21). He is currently a Visiting Professor CEBS, Mumbai. (emailid: physics.sutra@gmail.com)}
\vskip 0.5cm

\textit{Aneesh Kumar is Standard XII student at the Dhirubhai Ambani International School, Mumbai, who is keenly interested in Mathematics and Physics. (emailid: aneesh.kumar11235@gmail.com)}

\end{document}